\documentclass[12pt]{article}%
\usepackage[intlimits]{amsmath}
\usepackage{amssymb}
\usepackage[T1]{fontenc}
\usepackage[sc]{mathpazo}
\usepackage{color}
\usepackage[colorlinks]{hyperref}
\usepackage{amsfonts}
\usepackage{graphicx}%
\definecolor {refcol}{RGB}{40,0,255}
\hypersetup{colorlinks=true,allcolors=refcol}
\makeatletter
\newfont{\footsc}{cmcsc10 at 8truept}
\newfont{\footbf}{cmbx10 at 8truept}
\newfont{\footrm}{cmr10 at 10truept}
\newtheorem{theorem}{Theorem}

\newtheorem{corollary}[theorem]{Corollary}

\newtheorem{definition}[theorem]{Definition}

\newtheorem{problem}[theorem]{Problem}

\newenvironment{proof}[1][Proof]{\noindent{\textbf {#1}  }}  {\hfill$\Box$\bigskip}
\begin{document}

\title{\textbf{The }$p$-\textbf{norm of hypermatrices with symmetries}}
\author{V. Nikiforov\thanks{Department of Mathematical Sciences, University of
Memphis, Memphis, TN 38152, USA. Email: \textit{vnikifrv@memphis.edu}}}
\date{}
\maketitle

\begin{abstract}
The $p$-norm of $r$-matrices generalizes the $2$-norm of $2$-matrices. It is
shown that if a nonnegative $r$-matrix is symmetric with respect to two
indices $j$ and $k$, then the $p$-norm is attained for some set of vectors
such that the $i$th and the $j$th vectors are identical. It follows that the
$p$-spectral radius of a symmetric nonnegative $r$-matrix is equal to its
$p$-norm for any $p\geq2$.\medskip

\textbf{Keywords: }$p$-norm\textit{; hypermatrix; }$p$-\textit{spectral
radius; nonnegative hypermatrix; hypergraph. }

\textbf{AMS classification: }\textit{05C50, 05C65, 15A18, 15A42, 15A60,
15A69.}

\end{abstract}

\section{Introduction and main results}

In this note we study the $p$-norm of hypermatrices with partial symmetries.
To put our results in a more familiar context, we start with ordinary matrices.

Let $A=\left[  a_{i,j}\right]  $ be a real $m\times n$ matrix and let
\[
L_{A}\left(  \mathbf{x},\mathbf{y}\right)  :=\sum_{i,j}a_{i,j}x_{i}y_{j}%
\]
for any $\mathbf{x}=\left(  x_{1},\ldots,x_{m}\right)  \in\mathbb{R}^{m}$ and
$\mathbf{y}=\left(  y_{1},\ldots,y_{n}\right)  \in\mathbb{R}^{n}.$

Now, given $p\geq1,$ define the $p$-norm $\left\Vert A\right\Vert _{p}$ of $A$
as
\[
\left\Vert A\right\Vert _{p}:=\max_{\left\vert \mathbf{x}\right\vert
_{p}=1,\text{ }\left\vert \mathbf{y}\right\vert _{p}=1\text{ }}\text{ }%
|L_{A}\left(  \mathbf{x},\mathbf{y}\right)  |.
\]
where $\left\vert \cdot\right\vert _{p}$ stands for the $l^{p}$ vector norm.

It is known that if $A$ is symmetric, then
\[
\left\Vert A\right\Vert _{2}=\max_{\left\vert \mathbf{x}\right\vert
_{2}=1\text{ }}\text{ }|L_{A}\left(  \mathbf{x},\mathbf{x}\right)  |,
\]
but this equality turns out to be truly exceptional and fails for any $p\neq2$
and some appropriate matrix $A$. Our first theorem gives a condition on $A$
that is sufficient to preserve the equality for all $p\geq2.$

\begin{theorem}
\label{th2p}If $A$ is a symmetric nonnegative matrix and $p\geq2,$ then%
\[
\left\Vert A\right\Vert _{p}=\max_{\left\vert \mathbf{x}\right\vert _{p}%
=1}\text{ }L_{A}\left(  \mathbf{x},\mathbf{x}\right)  .
\]

\end{theorem}

The main goal of this note is to extend Theorem \ref{th2p} to $r$-matrices,
which we introduce next.

Let $r\geq2,$ and let $n_{1},\ldots,n_{r}$ be positive integers. An
$r$\emph{-matrix} of order $n_{1}\times\cdots\times n_{r}$ is a real function
defined on the Cartesian product $\left[  n_{1}\right]  \times\cdots
\times\left[  n_{r}\right]  .$

Thus, hereafter, \emph{matrix }means an $r$-matrix with unspecified $r$, and
ordinary matrices are referred to as $2$-matrices.

As usual, matrices are denoted by capital letters, and their values are
denoted by the corresponding lowercase letter with the variables listed as
subscripts, i.e., if $A$ is an $r$-matrix of order $n_{1}\times\cdots\times
n_{r}$, we write $a_{i_{1},\ldots,i_{r}}$ for $A\left(  i_{1},\ldots
,i_{r}\right)  $ whenever $i_{1}\in\left[  n_{1}\right]  ,\ldots,i_{r}%
\in\left[  n_{r}\right]  .$

The \emph{linear form}\textbf{ }of an $r$-matrix $A$ of order $n_{1}%
\times\cdots\times n_{r}$ is a function
\[
L_{A}:\mathbb{R}^{n_{1}}\times\cdots\times\mathbb{R}^{n_{r}}\rightarrow
\mathbb{R}%
\]
defined for any $r$ vectors \
\[
\mathbf{x}^{\left(  1\right)  }=(x_{1}^{\left(  1\right)  },\ldots,x_{n_{1}%
}^{\left(  1\right)  })\in\mathbb{R}^{n_{1}},\text{ }\ldots\text{ }%
,\mathbf{x}^{\left(  r\right)  }=(x_{1}^{\left(  r\right)  },\ldots,x_{n_{r}%
}^{\left(  r\right)  })\in\mathbb{R}^{n_{r}}%
\]
as
\[
L_{A}(\mathbf{x}^{\left(  1\right)  },\ldots,\mathbf{x}^{\left(  r\right)
}):=\sum_{i_{1}\in\left[  n_{1}\right]  ,\text{ }\ldots\text{ },i_{r}%
\in\left[  n_{r}\right]  }a_{i_{1},\ldots,i_{r}}x_{i_{1}}^{\left(  1\right)
}\text{ }\cdots\text{ }x_{i_{r}}^{\left(  r\right)  }.
\]
Now, the $p$\emph{-norm}\footnote{The idea of $p$-norm for hypermatrices comes
from Hardy, Littlewood and Polya \cite{HLP88}. In the present form it was
introduced for integral $p$ by Lek-Heng Lim \cite{Lim05}.}\textbf{ }of $A$ is
defined as%
\[
\left\Vert A\right\Vert _{p}:=\max_{|\mathbf{x}^{\left(  1\right)  }%
|_{p}=1,\text{ }\ldots\text{ },|\mathbf{x}^{\left(  r\right)  }|_{p}=1}\text{
}|L_{A}(\mathbf{x}^{\left(  1\right)  },\ldots,\mathbf{x}^{\left(  r\right)
})|.
\]

So far we have extended all of the setup of Theorem \ref{th2p} to
$r$-matrices, except for the symmetry property, which is less straightforward
for $r\geq3$.

\begin{definition}
\textbf{ }An $r$-matrix $A=\left[  a_{i_{1},\ldots,i_{r}}\right]  $ of order
$n_{1}\times\cdots\times n_{r}$ is called $\left(  j,k\right)  $%
\emph{-symmetric }if $i\neq k,$ $n_{j}=n_{k},$ and $A$ is invariant under the
swap of $i_{j}$ and $i_{k}$.

If $A$ is $\left(  j,k\right)  $-symmetric for every $1\leq j<k\leq r,$ then
it is called \emph{symmetric}. If an $r$-matrix $A$ of order $n_{1}%
\times\cdots\times n_{r}$ is symmetric, then $n_{1}=\cdots=n_{r},$ and $n_{1}$
is called its order.
\end{definition}

We are ready now to state our main results.

\begin{theorem}
\label{thr2}If $A$ is a $\left(  j,k\right)  $-symmetric $r$-matrix, then%
\[
\left\Vert A\right\Vert _{2}=\max_{|\mathbf{x}^{\left(  1\right)  }%
|_{2}=1,\text{ }\ldots\text{ },|\mathbf{x}^{\left(  r\right)  }|_{2}=1,\text{
\ }\mathbf{x}^{\left(  j\right)  }=\mathbf{x}^{\left(  k\right)  }}\text{
}|L_{A}\left(  \mathbf{x}^{\left(  1\right)  },\ldots,\mathbf{x}^{\left(
r\right)  }\right)  |.
\]

\end{theorem}

Note that statements similar to Theorem \ref{thr2} have been studied for
almost nine decades by now, mostly in abstract normed spaces. In particular,
motivated by problem 73 of Mazur and Orlicz in \cite{Sco81}, Banach
\cite{Ban38} proved a general result that implies Theorem \ref{thr2} for
symmetric matrices\footnote{For newer proofs of Banach's result, see
\cite{Fri13} and \cite{PST07}, and for further results, see \cite{Din99},
\cite{PST07} and their references.}. Our proof is much simpler, almost an
observation, because we use tools which may not available in general normed spaces.

If $A$ is nonnegative, Theorem \ref{thr2} can be extended similarly to Theorem
\ref{th2p}.

\begin{theorem}
\label{thrp}If $A$ is a $\left(  j,k\right)  $-symmetric nonnegative
$r$-matrix and $p\geq2,$ then%
\[
\left\Vert A\right\Vert _{p}=\max_{|\mathbf{x}^{\left(  1\right)  }%
|_{p}=1,\text{ }\ldots\text{ },|\mathbf{x}^{\left(  r\right)  }|_{p}=1,\text{
}\mathbf{x}^{\left(  j\right)  }=\mathbf{x}^{\left(  k\right)  }}\text{ }%
L_{A}\left(  \mathbf{x}^{\left(  1\right)  },\ldots,\mathbf{x}^{\left(
r\right)  }\right)  .
\]

\end{theorem}

The nonnegativity of $A$ in Theorem \ref{thrp} is far from necessary. Indeed
let $A$ be a $\left(  j,k\right)  $-symmetric nonnegative $r$-matrix of order
$n_{1}\times\cdots\times n_{r}$ and set $n=n_{j}=n_{k}$. Let $\mathbf{s}%
=\left(  s_{1},\ldots,s_{n}\right)  $ be a $\pm1$ vector. Define an $r$-matrix
$B$ of order $n_{1}\times\cdots\times n_{r}$ by letting
\[
b_{i_{1},\ldots,i_{r}}=a_{i_{1},\ldots,i_{r}}s_{i_{j}}s_{i_{k}}.
\]
Clearly $B$ may have both positive and negative entries, but it is easy to see
that
\[
\left\Vert B\right\Vert _{p}=\max_{|\mathbf{x}^{\left(  1\right)  }%
|_{p}=1,\text{ }\ldots\text{ },|\mathbf{x}^{\left(  r\right)  }|_{p}=1,\text{
}\mathbf{x}^{\left(  j\right)  }=\mathbf{x}^{\left(  k\right)  }}\text{ }%
L_{B}\left(  \mathbf{x}^{\left(  1\right)  },\ldots,\mathbf{x}^{\left(
r\right)  }\right)  .
\]

Another simple consequence of Theorem \ref{thrp} is the following corollary,
which was proved for $p\geq r$ in \cite{Nik16}, by rather involved methods..

\begin{corollary}
\label{corp}If $A$ is a symmetric nonnegative $r$-matrix and $p\geq2,$ then%
\[
\left\Vert A\right\Vert _{p}=\max_{|\mathbf{x}|_{p}=1}\text{ }L_{A}\left(
\mathbf{x},\ldots,\mathbf{x}\right)  .
\]

\end{corollary}

For symmetric $r$-matrices and $p\geq1$, the value $\max_{|\mathbf{x}|_{p}=1}$
$|L_{A}\left(  \mathbf{x},\ldots,\mathbf{x}\right)  |$ is known as the
$p$\emph{-spectral radius} of $A$ and is denoted by $\rho^{\left(  p\right)
}\left(  A\right)  .$ The $p$-spectral radius\footnote{The idea of the
$p$-spectral radius can be traced back to Lusternik and Schnirelman
\cite{LuSh30}, later revived by Friedman and Wigderson \cite{FrWi95}. For
hypergraphs it was introduced by Keevash, Lenz, and Mubayi \cite{KLM13}.} has
been studied in some detail, particularly for the adjacency matrix of uniform
hypergraphs (see \cite{Nik16} and its references). Corollary \ref{corp} can be
used to obtain various lower bounds on $\rho^{\left(  p\right)  }\left(
A\right)  ,$ by choosing vectors $\mathbf{x}^{\left(  1\right)  }%
,\ldots,\mathbf{x}^{\left(  r\right)  }$ with $|\mathbf{x}^{\left(  1\right)
}|_{p}=1,$ $\ldots$ $,$ $|\mathbf{x}^{\left(  r\right)  }|_{p}=1$ and
calculating $L_{A}\left(  \mathbf{x}^{\left(  1\right)  },\ldots
,\mathbf{x}^{\left(  r\right)  }\right)  .$ We give an illustration next.

Suppose that $A$ is an $r$-matrix of order $n_{1}\times\cdots\times n_{r}$ and
for every $i\in\left[  n_{1}\right]  ,$ set%
\[
S_{i}=\sum_{i_{2}\in\left[  n_{2}\right]  ,\text{ }\ldots\text{ },i_{r}%
\in\left[  n_{r}\right]  }a_{i,i_{2}\ldots,i_{r}}.
\]
Apparently the values $S_{1},\ldots,S_{n_{1}}$ generalize the row-sums of
$2$-matrices, so we call them the \emph{slice-sums} of $A.$

Combining Theorem 21 of \ \cite{Nik16} with Theorem \ref{thrp}, one comes up
with the following useful lower bound on the $p$-spectral radius:

\begin{corollary}
Let $A$ be a symmetric nonnegative $r$-matrix of order $n$ with slice-sums
$S_{1},\ldots,S_{n}.$ If $p\geq2,$ then%
\[
\rho^{\left(  p\right)  }\left(  A\right)  \geq n^{1-r/p}\left(  \frac{1}%
{n}\left(  S_{1}^{p/\left(  p-1\right)  }+\cdots+S_{n}^{p/\left(  p-1\right)
}\right)  \right)  ^{\left(  p-1\right)  /p}.
\]

\end{corollary}

When restated for $r$-uniform graphs (see, e.g., \cite{Nik16} for the basics),
the corollary reads as:

\begin{corollary}
\label{cogr}Let $G$ be an $r$-uniform graph of order $n$ with degrees
$d_{1},\ldots,d_{n}.$ If $p\geq2,$ then%
\[
\rho^{\left(  p\right)  }\left(  G\right)  \geq\left(  r-1\right)
!n^{1-r/p}\left(  \frac{1}{n}\left(  d_{1}^{p/\left(  p-1\right)  }%
+\cdots+d_{n}^{p/\left(  p-1\right)  }\right)  \right)  ^{\left(  p-1\right)
/p}.
\]

\end{corollary}

Let us note that the case $p=r$ of this corollary has been proved in
\cite{LKS16}, and the case $p\geq r$ has been proved in \cite{Nik16}.

The remaining part of the note is split into two sections: in Section
\ref{pfs}, we present the proofs of Theorems \ref{th2p}, \ref{thr2},
\ref{thrp} and in Section \ref{pros}, we state two open problems.

\section{\label{pfs}Proofs of Theorems \ref{th2p}, \ref{thr2}, and \ref{thrp}}

\begin{proof}
[\textbf{Proof of Theorem \ref{th2p}}]Let $A$ be a symmetric nonnegative
matrix of order $n,$ and
\[
\left\Vert A\right\Vert _{p}=\left\vert L_{A}\left(  \mathbf{x},\mathbf{y}%
\right)  \right\vert
\]
where $\left\vert \mathbf{x}\right\vert _{p}=\left\vert \mathbf{y}\right\vert
_{p}=1.$ We assume that $\mathbf{x}$ and $\mathbf{y}$ are nonnegative because
\[
|\sum_{i,j}a_{i,j}x_{i}y_{j}|\text{ }\leq\sum_{i,j}a_{i,j}\left\vert
x_{i}\right\vert \left\vert y_{j}\right\vert .
\]
and the $l^{p\text{ }}$norm of both $\left(  \left\vert x_{1}\right\vert
,\ldots,\left\vert x_{n}\right\vert \right)  $ and $\left(  \left\vert
y_{1}\right\vert ,\ldots,\left\vert y_{n}\right\vert \right)  $ is $1.$

We see that $\mathbf{x}$ and $\mathbf{y}$ maximize $L_{A}\left(
\mathbf{x},\mathbf{y}\right)  $ under the constraints $x_{1}^{p}+\cdots
+x_{n}^{p}=1$ and $y_{1}^{p}+\cdots\cdots+y_{n}^{p}=1.$ Using Lagrange's
multipliers, it follows that there exists $\lambda$ such that
\begin{equation}
\lambda x_{i}^{p-1}=a_{i,1}y_{1}+\cdots+a_{i,n}y_{n},\text{ }i=1,\ldots,n
\label{eqx}%
\end{equation}
and%
\begin{equation}
\lambda y_{i}^{p-1}=a_{i,1}x_{1}+\cdots+a_{i,n}x_{n},\text{ }i=1,\ldots,n
\label{eqy}%
\end{equation}
It is easy to see that $\lambda=\left\Vert A\right\Vert _{p}:$ indeed,
multiplying the $i$th equation (\ref{eqx}) by $x_{i}$ and adding the results,
we get%
\[
\lambda=\lambda\left(  x_{1}^{p}+\cdots+x_{n}^{p}\right)  =\sum_{i,j}%
a_{i,j}x_{i}y_{j}=\left\Vert A\right\Vert _{p}.
\]

Now, let
\[
z_{i}:=\left(  \frac{x_{i}^{p}+y_{i}^{p}}{2}\right)  ^{1/p},\text{
\ \ \ }i=1,\ldots,n.
\]
Multiplying the $i$th equation (\ref{eqx}) by $x_{i}$ and the $i$th equation
(\ref{eqy}) by $y_{i}$ and adding the results, we get
\[
\sum_{i=1}^{n}\lambda x_{i}^{p}+\lambda y_{i}^{p}=\sum_{i=1}^{n}\sum_{j=1}%
^{n}a_{i,j}\left(  x_{i}y_{j}+y_{i}x_{j}\right)  .
\]
On the other hand, using the Cauchy-Schwarz inequality and the Power Mean
inequality, we find that
\[
\frac{x_{i}y_{j}+y_{i}x_{j}}{2}\leq\left(  \frac{x_{i}^{2}+y_{i}^{2}}%
{2}\right)  ^{1/2}\left(  \frac{x_{j}^{2}+y_{j}^{2}}{2}\right)  ^{1/2}%
\leq\left(  \frac{x_{i}^{p}+y_{i}^{p}}{2}\right)  ^{1/p}\left(  \frac
{x_{j}^{p}+y_{j}^{p}}{2}\right)  ^{1/p}=z_{i}z_{j}.
\]
Therefore,
\[
\left\Vert A\right\Vert _{p}=\lambda=\sum_{i=1}^{n}\lambda\frac{x_{i}%
^{p}+y_{i}^{p}}{2}\leq\sum_{i=1}^{n}\sum_{j=1}^{n}a_{i,j}z_{i}z_{j}%
=L_{A}\left(  \mathbf{z},\mathbf{z}\right)  .
\]
Since $\left\vert \mathbf{z}\right\vert _{p}=1,$ the proof is completed.
\end{proof}

\bigskip

\begin{proof}
[\textbf{Proof of Theorem \ref{thr2}}]Let $A$ be a $\left(  j,k\right)
$-symmetric $r$-matrix of order $n_{1}\times\cdots\times n_{r}$. For
convenience, let us reindex the variables so that $j=r-1,$ $k=r,$ and let
$n=n_{r-1}=n_{r}.$

Suppose that $\mathbf{x}^{\left(  1\right)  }\in\mathbb{R}^{n_{1}},$ $\ldots$
$,\mathbf{x}^{\left(  r\right)  }\in\mathbb{R}^{n_{r}}$ are vectors with
$|\mathbf{x}^{\left(  1\right)  }|_{2}=1,$ $\ldots$ $,|\mathbf{x}^{\left(
r\right)  }|_{2}=1$ such that
\[
\left\Vert A\right\Vert _{2}=|L_{A}(\mathbf{x}^{\left(  1\right)  }%
,\ldots,\mathbf{x}^{\left(  r\right)  })|.
\]
To finish the proof we have to show that $\mathbf{x}^{\left(  r-1\right)  }$
and $\mathbf{x}^{\left(  r\right)  }$ may be chosen equal. To this end, define
a square $2$-matrix $B$ of order $n$ by%
\[
b_{s,t}=\sum_{i_{1}\in\left[  n_{1}\right]  ,\text{ }\ldots\text{ },i_{r-2}%
\in\left[  n_{r-2}\right]  }a_{i_{1},\ldots i_{r-2},s,t}x_{i_{1}}^{\left(
1\right)  }\text{ }\cdots\text{ }x_{i_{r-2}}^{\left(  r-2\right)  }.
\]
and note that $B$ is symmetric since $A$ is $\left(  r-1,r\right)  $-symmetric.

Next, for convenience, set $\mathbf{x}=\mathbf{x}^{\left(  r-1\right)  }$ and
$\mathbf{y}=\mathbf{x}^{\left(  r\right)  }.$ Obviously%
\[
L_{A}(\mathbf{x}^{\left(  1\right)  },\ldots,\mathbf{x}^{\left(  r\right)
})=\sum_{i=1}^{n}\sum_{j=1}^{n}b_{i,j}x_{i}y_{j}=L_{B}(\mathbf{x},\mathbf{y}).
\]
Clearly $\mathbf{x}$ and $\mathbf{y}$ maximize $L_{B}(\mathbf{x},\mathbf{y})$
subject to $x_{1}^{2}+\cdots+x_{n}^{2}=1$ and $y_{1}^{2}+\cdots+y_{n}^{2}=1.$
Therefore,%
\[
\left\Vert A\right\Vert _{2}=\left\Vert B\right\Vert _{2}.
\]
Moreover, using Lagrange's multipliers, it follows that there exists $\lambda$
such that
\begin{equation}
\lambda x_{i}=b_{i,1}y_{1}+\cdots+b_{i,n}y_{n},\text{ \ \ \ \ }i=1,\ldots
,n\label{eqx2}%
\end{equation}
and%
\begin{equation}
\lambda y_{i}=b_{i,1}x_{1}+\cdots+b_{i,n}x_{n},\text{ \ \ \ \ }i=1,\ldots
,n\label{eqy2}%
\end{equation}

As in the proof of Theorem \ref{th2p}, we see that $\lambda=\left\Vert
B\right\Vert _{2},$ and so $\lambda=\left\Vert A\right\Vert _{2}.$

Further, adding the $i$th equation (\ref{eqx2}) and the $i$th equation
(\ref{eqy2}), we get a new system of equations%
\[
\lambda\left(  x_{i}+y_{i}\right)  =a_{i,1}\left(  x_{1}+y_{1}\right)
+\cdots+a_{i,n}\left(  x_{n}+y_{n}\right)  ,\text{ \ \ \ }i=1,\ldots,n
\]
which implies that $\lambda$ is an eigenvalue of $B,$ unless $\mathbf{x}%
+\mathbf{y}=0.$ In the latter case, we immediately see that
\[
\left\Vert A\right\Vert _{2}=|L_{A}(\mathbf{x}^{\left(  1\right)  }%
,\ldots,\mathbf{x}^{\left(  r-2\right)  },\mathbf{x},\mathbf{x})|,
\]
completing the proof. Thus, we may assume that $\lambda$ is an eigenvalue of
$B.$ The Rayleigh-Ritz theorem implies that
\[
\lambda\leq\max_{\left\vert \mathbf{z}\right\vert _{2}=1}\text{ }|L_{B}\left(
\mathbf{z},\mathbf{z}\right)  |.
\]
Hence,%
\[
\left\Vert A\right\Vert _{2}\leq\max_{\left\vert \mathbf{z}\right\vert _{2}%
=1}\text{ }|L_{B}\left(  \mathbf{z},\mathbf{z}\right)  |\text{ }%
=\max_{\left\vert \mathbf{z}\right\vert _{2}=1}\text{ }|L_{A}(\mathbf{x}%
^{\left(  1\right)  },\ldots,\mathbf{x}^{\left(  r-2\right)  },\mathbf{z}%
,\mathbf{z})|,
\]
completing the proof.
\end{proof}

\bigskip

\begin{proof}
[\textbf{Proof of Theorem \ref{thrp}}]Our proof is a combination of the proofs
of Theorems \ref{th2p} and \ref{thr2}.

Let $A$ be a $\left(  j,k\right)  $-symmetric nonnegative $r$-matrix of order
$n_{1}\times\cdots\times n_{r}$. For convenience, let us reindex the variables
so that $j=r-1,$ $k=r,$ and let $n=n_{r-1}=n_{r}.$

Suppose that $\mathbf{x}^{\left(  1\right)  }\in\mathbb{R}^{n_{1}},$ $\ldots$
$,\mathbf{x}^{\left(  r\right)  }\in\mathbb{R}^{n_{r}}$ are vectors with
$|\mathbf{x}^{\left(  1\right)  }|_{p}=1,$ $\ldots$ $,|\mathbf{x}^{\left(
r\right)  }|_{p}=1$ such that
\[
\left\Vert A\right\Vert _{p}=|L_{A}(\mathbf{x}^{\left(  1\right)  }%
,\ldots,\mathbf{x}^{\left(  r\right)  })|.
\]
Note that $\mathbf{x}^{\left(  1\right)  },$ $\ldots$ $,\mathbf{x}^{\left(
r\right)  }$ can be taken nonnegative in view of
\[
|L_{A}(\mathbf{x}^{\left(  1\right)  },\ldots,\mathbf{x}^{\left(  r\right)
})|\text{ }\leq L_{A}(|\mathbf{x}^{\left(  1\right)  }|,\ldots,|\mathbf{x}%
^{\left(  r\right)  }|)\text{ }%
\]
To finish the proof we have to show that $\mathbf{x}^{\left(  r-1\right)  }$
and $\mathbf{x}^{\left(  r\right)  }$ can be chosen equal. To this end, define
a square $2$-matrix $B$ of order $n$ by%
\[
b_{s,t}=\sum_{i_{1}\in\left[  n_{1}\right]  ,\text{ }\ldots\text{ },i_{r-2}%
\in\left[  n_{r-2}\right]  }a_{i_{1},\ldots i_{r-2},s,t}x_{i_{1}}^{\left(
1\right)  }\cdots x_{i_{r-2}}^{\left(  r-2\right)  }%
\]
and note that $B$ is symmetric since $A$ is $\left(  r-1,r\right)
$-symmetric. Moreover, $B$ is nonnegative.

Next, for convenience, set $\mathbf{x}=\mathbf{x}^{\left(  r-1\right)  }$ and
$\mathbf{y}=\mathbf{x}^{\left(  r\right)  }.$ Obviously%
\[
L_{A}(\mathbf{x}^{\left(  1\right)  },\ldots,\mathbf{x}^{\left(  r\right)
})=\sum_{i=1}^{n}\sum_{j=1}^{n}b_{i,j}x_{i}y_{j}=L_{B}(\mathbf{x},\mathbf{y}).
\]
Clearly $\mathbf{x}$ and $\mathbf{y}$ maximize $L_{B}(\mathbf{x},\mathbf{y})$
subject to $x_{1}^{p}+\cdots+x_{n}^{p}=1$ and $y_{1}^{p}+\cdots+y_{n}^{p}=1.$
Therefore,%
\[
\left\Vert A\right\Vert _{p}=\left\Vert B\right\Vert _{p}.
\]

In view of Theorem \ref{th2p}, there exists $\mathbf{z}\in\mathbb{R}^{n}$ with
$\left\vert \mathbf{z}\right\vert _{p}=1$ such that
\[
\left\Vert B\right\Vert _{p}=L_{B}(\mathbf{z},\mathbf{z}).
\]
Hence,%
\[
L_{A}(\mathbf{x}^{\left(  1\right)  },\ldots,\mathbf{x}^{\left(  r-2\right)
},\mathbf{z},\mathbf{z})=L_{B}(\mathbf{z},\mathbf{z})=\left\Vert B\right\Vert
_{p}=\left\Vert A\right\Vert _{p},
\]
completing the proof.
\end{proof}

\section{\label{pros}Two open problems}

Corollary \ref{corp} could be very useful because $\left\Vert A\right\Vert
_{p}$ is usually easier to evaluate or estimate than $\rho^{\left(  p\right)
}\left(  A\right)  .$ Thus, it is desirable to extend Corollary \ref{corp} to
matrices that are essentially distinct from nonnegative matrices. We naturally
arrive at the following problems:

\begin{problem}
Characterize all symmetric $r$-matrices $A$ such that
\[
\rho^{\left(  p\right)  }\left(  A\right)  =\left\Vert A\right\Vert _{p}%
\]
for all sufficiently large $p$.
\end{problem}

\begin{problem}
Characterize all symmetric $r$-matrices $A$ such that
\[
\rho^{\left(  p\right)  }\left(  A\right)  =\left\Vert A\right\Vert _{p}%
\]
for all $p\in\left[  1,2\right)  .$
\end{problem}

Unfortunately, the above problems seem hopeless at present. Probably there are
some chances for solving either of them for $r=2.$

\end{document}